%

\documentstyle{amsppt}
\loadbold
\magnification=1200
\NoRunningHeads
\def\vDash{\mathrel{\|}\joinrel\mathrel{-}}
\topmatter
\title
Remark on the failure of \\ Martin's Axiom
\endtitle
\author
Avner Landver
\endauthor
\address
Department of Mathematics,
The University of Kansas,
Lawrence, KS 66045
\endaddress
\email
landver\@kuhub.cc.ukans.edu
\endemail
\endtopmatter

\document
\baselineskip 18pt
\parskip 6pt

Let $m$ be the least cardinal $\theta$ such that MA$_\theta$ fails.
The only known model for ``$m$ is singular" was constructed by Kunen [K1].
In Kunen's model $cof(m)=\omega_1$.
It is unknown whether $``\omega_1 < cof(m) < m"$ is consistent.
The purpose of this paper is to present a proof of Kunen's result and  
to identify the difficulties of generalizing this result to an arbitrary
uncountable cofinality. The following material
is based on [K1]. We would like to thank K. Kunen and S. Shelah [S] for their 
input.  

\subhead \S0  Definitions and some facts \endsubhead

For undefined terminology consult [K].
Let  $\theta$ be a fixed singular cardinal with uncountable cofinality.
Let $\kappa=cof(\theta)$, and fix $\{\theta_\alpha : \alpha<\kappa\}$
an increasing sequence of cardinals converging to $\theta$ with 
$\theta_0 > \kappa$.

Let ${\Cal A}=\{a_\xi : \xi\in\theta\}\subseteq\wp(\omega)$ 
and ${\Cal B}=\{b_\xi : \xi\in\theta\}\subseteq\wp(\omega)$. 
$({\Cal A},{\Cal B})$ is a {\bf $(\theta,\theta)$ pair} if 
for every $ \xi,\eta\in\theta$, $|a_\xi \cap b_\eta|<\aleph_0$.
The pair $({\Cal A},{\Cal B})$ is {\bf disjoint} if for every $\xi \in \theta$,
$a_{\xi} \cap b_{\xi} = \emptyset$.
The pair $({\Cal A},{\Cal B})$ is {\bf locally split} if for every ${\Cal X}
\subset {\Cal A}$, and every ${\Cal Y}\subset{\Cal B}$, such that
$| {\Cal X} \cup {\Cal Y}| \le \aleph_1$, the pair $({\Cal X},{\Cal Y})$ 
{\bf splits} (i.e. $(\exists c \subseteq \omega) [(\forall a \in {\Cal X})
|a \setminus c|< \aleph_0$ and $(\forall b \in {\Cal Y}) |c \cap b| < 
\aleph_0]$). 
Notice that if $({\Cal A},{\Cal B})$ is locally split in V, then 
$({\Cal A},{\Cal B})$ is locally split in every c.c.c.\ extension of V.

For a disjoint pair $({\Cal A},{\Cal B})$ we define the partial order 
$$
S({\Cal A},{\Cal B})=
\{f:\theta \rightarrow \omega : |f| < \aleph_0 \,\,\text{and}\,\,
(\forall \xi,\eta\in dom(f)) \, (a_\xi - f(\xi))\cap (b_\eta - f(\eta))=
\emptyset \}.
$$
$S({\Cal A},{\Cal B})$ is partially ordered by inclusion.
It is well known that if $({\Cal A},{\Cal B})$ is locally split, then
$S({\Cal A},{\Cal B})$ is c.c.c., and therefore, if in addition MA$_\theta$ 
holds, then 
$({\Cal A},{\Cal B})$ can be split (to split $({\Cal A},{\Cal B})$ use the set
$c=\bigcup \{a_\xi - f(\xi) : f \in G, \xi \in dom(f)\}$
Where $G$ is a sufficiently generic filter on $S({\Cal A},{\Cal B})$).

If $x \in [\theta ]^{<\aleph_0}$, then we define 
$$
a(x)=\bigcap_{\xi\in x} a_\xi.
$$ 
$a( \emptyset )= \omega$, and $b(x)$ is defined analogously. 

The pair $({\Cal A},{\Cal B})$ is {\bf nice} if 
for every $X=\{x_i : i\in\theta\}$, a family of
disjoint finite subsets of $\theta$, there are $i,j \in \theta$
such that $a(x_i)\cap b(x_j)\neq \emptyset$. It is not hard to see that a nice 
disjoint $(\theta , \theta)$ pair is also a {\bf gap} (i.e. can not be split).  
Therefore, the existence
of a nice, locally split, disjoint $(\theta , \theta)$ pair contradicts 
MA$_\theta$. ``Nice" is a variation of Kunen's ``strong gap" [K1].  

To add a nice, disjoint $(\theta , \theta)$ gap that is locally split, one uses 
the partial order ${\Cal P}$ defined by: 
$$
p \in {\Cal P} \iff p=\langle n_p, s_p, a_p, b_p, Z_p, 
\langle c_p^z : z\in Z_p \rangle \rangle
$$
where $n_p \in \omega$, $s_p \in [\theta]^{<\aleph_0}$, $a_p$ \& $b_p \subseteq 
n_p \times s_p$, $a_p \cap b_p = \emptyset$, $Z_p \in [[\theta]^{\le\aleph_1}
]^{<\aleph_0}$, and $(\forall z \in Z_p)\, c_p^z \subseteq n_p$. 
The ordering on $\Cal P$ is defined by putting 
$q \leq p$ if and only if the following five conditions hold:
$$
\align
&{\tenrm(1)}\,s_q \supseteq s_p \wedge n_q \geq n_p \wedge Z_q \supseteq Z_p.\\ 
&{\tenrm(2)}\,a_q \cap (n_p \times s_p) = a_p \,\,\text{and}\,\, 
b_q \cap (n_p \times s_p) = b_p.\\
&{\tenrm(3)}\,(\forall z\in Z_p)\, c_q^z \cap n_p = c_p^z.\\
&{\tenrm(4)}\,(\forall \xi , \eta \in s_p)(\forall l \in [n_p,n_q))\, 
[\langle l,\xi \rangle \in 
a_q \rightarrow \langle l,\eta \rangle \not\in b_q].\\
&{\tenrm(5)}\,(\forall z\in Z_p)(\forall \xi \in z\cap s_p)(\forall l\in [n_p,n_q)
\,[(\langle l,\xi \rangle \in a_q \rightarrow l \in c_q^z)\,\,\wedge \,\,
(\langle l,\xi \rangle \in b_q \rightarrow  l \not\in c_q^z)].
\endalign
$$
${\Cal P}$ is c.c.c., and forcing with ${\Cal P}$ adds the disjoint pair 
${\Cal A}=\{a_\xi : \xi\in\theta\}, {\Cal B}=\{b_\xi : \xi\in\theta\}$,
where $a_\xi = \{l\in \omega : (\exists p\in G)\, \langle l,\xi \rangle
\in a_p \}$ ($b_\xi$ is defined similarly), and $G$ is the generic filter.
To see that $({\Cal A},{\Cal B})$ is locally split it is enough to show that
for every $z \in [\theta]^{\aleph_1} \cap \text{V}$, the pair
$(\{a_\xi : \xi\in z\}, \{b_\xi : \xi\in z\})$ splits (this is enough 
because ${\Cal P}$ is c.c.c., and therefore every subset of $\theta$ with
cardinality $\aleph_1$ in the extension is contained in such a subset from V).
But this pair get split by $\bigcup \{c_p^z : p \in G \,\,
\text{and}\,\, z \in Z_p \}$.  Finally let us show that

\proclaim {Fact 0}
The pair $({\Cal A},{\Cal B})$ is nice.
\endproclaim

\demo {Proof}
Assume that $p \in {\Cal P}$ is such that 
$$
p \Vdash ``X=\{x_i : i\in\theta\}\,\,\text{are finite disjoint  
subsets of}\,\,\theta".
$$  
For every $i \in \theta$, let $p_i \leq p$ and $y_i \in [\theta]^{<\aleph_0}$
be such that $p_i \Vdash ``x_i = y_i"$.
Let $\kappa < \lambda < \theta$ be an arbitrary regular cardinal.
Let $A \in [\theta]^\lambda$ and $n \in \omega$ be such that:\newline
(a) $(\forall i \in A) n_{p_i}=n$. \newline
(b) $\{s_{p_i} : i \in A\}$ form a delta system with root $s$ and 
$(\forall i,j \in A) [a_{p_i}\cap (n \times s)=a_{p_j}\cap (n \times s)$
and $b_{p_i}\cap (n \times s)=b_{p_j}\cap (n \times s)]$. \newline
(c) $\{Z_{p_i} : i \in A\}$ form a delta system with root $Z$ and
$(\forall z \in Z) c_{p_i}^z = c_{p_j}^z$. \newline
(d) $(\forall i \in A) s_{p_i} \supset y_i$.

Notice that (a)-(c) imply that $\{p_i : i\in A\}$ are linked (i.e. pairwise
compatible) and therefore $\{y_i : i\in A\}$ are disjoint.
Let $\bar Z = \bigcup Z$, then  $|\bar Z| \leq \aleph_1$. Therefore we can
find $i \neq j \in A$ such that
$$
y_i \cap s = y_i \cap \bar Z = y_j \cap s = y_j \cap \bar Z = \emptyset.
$$
Finally define $q \leq p_i,p_j$ as follows: let $n_q = n+1$, 
let $s_q = s_{p_i} \cup s_{p_j}$, and $Z_q = Z_{p_i} \cup Z_{p_j}$.
We put
$$
\align
&\langle n,\xi \rangle \in a_q \iff \xi \in y_i \\
&\langle n,\xi \rangle \in b_q \iff \xi \in y_j.
\endalign
$$
We also make sure that for every $z \in (Z_{p_i} \setminus Z)$,
if $\xi \in z \cap y_i$, then $n \in c_q^z$,
and that for every $z \in (Z_{p_j} \setminus Z)$, 
if $\xi \in z \cap y_j$, then $n \notin c_q^z$.
Notice that this does not cause a contradiction since only $z \notin Z$ are 
involved. We conclude that $q \leq p$ and that
$$
q \Vdash ``n \in a(x_i) \cap b(x_j)".  \qed
$$
\enddemo

Given a disjoint 
$(\theta , \theta)$ pair $({\Cal A},{\Cal B})$ and 
$X= \{x_i : i \in \theta \}$ a family of disjoint finite subsets of $\theta$,
we make the following definition. 
We call $S=\{S_{\alpha} : \alpha < \kappa \}$ an ${\bold X}${\bf -sequence} 
if
the $S_{\alpha}$'s are disjoint subsets of $\theta$ and 
$(\forall \alpha < \kappa)\, [\,| S_{\alpha} |>\theta _{\alpha} $ and
$(\forall i, j \in S_{\alpha}) \, a(x_i) \cap b(x_j) = \emptyset \, ]$.

Next, for every $X= \{x_i : i \in \theta \}$ disjoint finite subsets of $\theta$,
and every 
$S=\{S_{\alpha} : \alpha < \kappa \}$ an $X$-sequence,
we define: 
$$
\align
{\Cal Q}(X,S) = \{&F \in [\kappa ]^{< \aleph_0} : 
(\forall \alpha \neq \beta \in F) (\exists i \in S_{\alpha})
(\exists j \in S_{\beta}) \\
&[a(x_i) \cap b(x_j) \neq \emptyset \,\vee \, 
b(x_i) \cap a(x_j) \neq \emptyset ]\}
\endalign
$$
and
$$
\align
{\Cal P}(X,S) = \{&F \in [\kappa ]^{< \aleph_0} : 
(\forall \alpha \neq \beta \in F) (\forall i \in S_{\alpha})
(\forall j \in S_{\beta}) \\
&[a(x_i) \cap b(x_j) = \emptyset \,\wedge \, 
b(x_i) \cap a(x_j) = \emptyset ]\}.
\endalign
$$
${\Cal Q}(X,S)$ and ${\Cal P}(X,S)$ are both partially ordered by inclusion.
Notice that the last three definitions depend on 
${\Cal A}=\{a_\xi : \xi\in\theta\}$ 
and ${\Cal B}=\{b_\xi : \xi\in\theta\}$.
${\Cal P}(X,S)$ is a typical ``dangerous" partial order (see definition 1),
and ${\Cal Q}(X,S)$ will be used in the proof of Kunen's result to ``kill"
dangerous partial orders.  Notice that ${\Cal P}(X,S) \times {\Cal Q}(X,S)$
is not $\kappa$.c.c\. 
(the set $\{\langle\{\alpha\},\{\alpha\}\rangle: \alpha < \kappa\}$ is an 
antichain of size $\kappa$).

\subhead \S1 A proof of Kunen's result \endsubhead

The first step is to show that the niceness of the pair is in fact a statement
concerning the $\kappa$.c.c\. of the various ${\Cal Q}(X,S)$'s.

\proclaim {Lemma 1} 
Assume that $({\Cal A}, {\Cal B})$ is a $(\theta ,\theta)$ pair. 
Then $({\Cal A}, {\Cal B})$ is nice if and only if for every 
$X= \{x_u : u \in \theta \}$, $Y= \{y_v : v \in \theta \}$
disjoint finite subsets of $\theta$
and every $S=\{S_{\alpha} : \alpha < \kappa \}$ an $X$-sequence, and
$T=\{T_{\alpha} : \alpha < \kappa \}$ a $Y$-sequence, the partial order
${\Cal Q}(X,S) \times {\Cal Q}(Y,T)$ is $\kappa$.c.c.
\endproclaim

\demo {Proof}
The reverse implication is easy to check.  Let us prove the direct implication.
Assume that $\{\langle K_\gamma , F_\gamma \rangle : \gamma < \kappa \}$
is an antichain in  ${\Cal Q}(X,S) \times {\Cal Q}(Y,T)$.  
We may assume that $\{ K_\gamma : \gamma < \kappa \}$ are disjoint and all
have size $n$, and that the $F_\gamma$'s are disjoint and all have size $p$ 
(use delta systems).

For $t \in [\theta]^{<\aleph_0}$, we call $t$ a 
${\bold \gamma}$-{\bf transversal} if
$|t|=n+p$, and $(\forall\alpha \in K_\gamma)\, |t \cap S_\alpha|=1$,
and $(\forall\beta \in F_\gamma)\, |t \cap T_\beta|=1$.
Choose $\{t_i : i \in \theta \}$ pairwise disjoint, with each $t_i$
being a $\gamma$-transversal for some $\gamma < \kappa$. 
Next, let 
$Z=\{z_i : i \in \theta \}$ be defined in the following way.
If $t_i$ is a $\gamma$-transversal, then
$$
z_i=(\cup \{x_u : u \in (t_i \cap S_\alpha) \wedge \alpha \in K_\gamma  \})\,
\cup \,(\cup \{y_v : v \in (t_i \cap T_\beta) \wedge \beta \in F_\gamma  \}).
$$
Since both $\{x_u : u \in \theta \}$ and 
$\{y_v : v \in \theta \}$ are pairwise disjoint, and the transversals
$\{t_i : i \in \theta \}$ are pairwise disjoint, we may assume that
$Z=\{z_i : i \in \theta \}$ are pairwise disjoint.

Finally, we show that the existence of $Z$ 
contradicts the niceness of $({\Cal A}, {\Cal B})$.
Let $i \neq j \in \theta$.

{\it Case} 1: $t_i , t_j$ are both $\gamma$-transversals.  Let
$\alpha \in K_\gamma$ (if $n=0$, then work with the $F_\gamma$'s).
Let $u \in (t_i \cap S_\alpha)$ and $w \in (t_j \cap S_\alpha)$
with $u \neq w$. 
Clearly, $a(x_u) \cap b(x_w) = b(x_u) \cap a(x_w) = \emptyset$
(because $S$ is an $X$-sequence).
But $x_u \subset z_i$ and  $x_w \subset z_j$
therefore $a(z_i) \cap b(z_j) = b(z_i) \cap a(z_j) = \emptyset$.
\newline

{\it Case} 2: $t_i$ is a $\gamma$-transversal, and 
$t_j$ is a $\delta$-transversal, and $\gamma \neq \delta$.
In this case $\langle K_\gamma , F_\gamma \rangle \perp 
\langle K_\delta , F_\delta \rangle$.  Assume w.l.o.g.\ that
$K_\gamma \perp K_\delta$.  This means that 
$$
(\exists \alpha \in K_\gamma)(\exists \beta \in K_\delta)
(\forall u \in S_\alpha)(\forall w \in S_\beta)\,
a(x_u) \cap b(x_w) = b(x_u) \cap a(x_w) = \emptyset.
$$
Now let $u \in (t_i \cap S_\alpha)$ and $ w \in (t_j \cap S_\beta)$.
By the above, $a(x_u) \cap b(x_w) = b(x_u) \cap a(x_w) = \emptyset$.
But $x_u \subset z_i$ and  $x_w \subset z_j$,
therefore $a(z_i) \cap b(z_j) = b(z_i) \cap a(z_j) = \emptyset$.
\qed
\enddemo

Similarly, it can be shown that $({\Cal A}, {\Cal B})$ is nice if and only if
every ${\Cal Q}(X,S)$ is $\kappa$.c.c.\ , and also if and only if the product of 
any finitely many partial orders of the form ${\Cal Q}(X,S)$ is $\kappa$.c.c. 

It is also true that if $\tilde S=\{S_{\alpha} : \alpha < \kappa \}$ is defined 
by $S_{\alpha}=[\theta_\alpha, \theta_\alpha^+)$, then: 
$({\Cal A}, {\Cal B})$ is nice if and only if for every 
$X= \{x_u : u \in \theta \}$ disjoint finite subsets of $\theta$,
if $\tilde S$ is an $X$-sequence, then
${\Cal Q}(X,\tilde S)$ is $\kappa$.c.c.

\definition {Definition 1}
Assume that $({\Cal A}, {\Cal B})$ is a nice pair.  The partial order ${\Cal R}$
is called {\bf dangerous} for $({\Cal A}, {\Cal B})$ if 
there exsits $r \in {\Cal R}$ such that 
$$
r\vDash_{\Cal R} ``({\Cal A}, {\Cal B})\text{ is
not nice} ". 
$$ 
\enddefinition

The following style of proof was motivated by [S].

\proclaim {Lemma 2} 
Assume that $({\Cal A}, {\Cal B})$ is a nice $(\theta ,\theta)$ pair, and 
${\Cal R}$ is a  $\kappa$.c.c.\ partial order with $|{\Cal R}| < \theta$.  Then 
${\Cal R}$ is dangerous for $({\Cal A}, {\Cal B})$
if and only if 
there exist 
$Y= \{y_i : i \in \theta \}$ disjoint finite subsets of $\theta$,
and $T=\{T_{\alpha} : \alpha < \kappa \}$ a $Y$-sequence,
such that ${\Cal R} \times {\Cal Q}(Y,T)$ is not $\kappa$.c.c.
\endproclaim

\demo{Proof}
($\Leftarrow$):
If ${\Cal R} \times {\Cal Q}(Y,T)$ is not $\kappa$.c.c., then
there exists $r \in {\Cal R}$ such that 
$$
r\vDash_{\Cal R} ``{\Cal Q}(Y,T)\text{ is 
not}\,\, \kappa\text{.c.c.\ }". 
$$
Therefore, by Lemma 1, ${\Cal R}$ is dangerous. \newline 
($\Rightarrow$):
Let $r \in {\Cal R}$ be such that
$$
r \vDash ``X=\{x_i : i \in \theta \}\,\text{ are disjoint finite
subsets of}\,\, \theta \text{ and }
(\forall i,j \in \theta) \,  a(x_i) \cap b(x_j) = \emptyset".
$$
$(\forall i \in \theta)$\,\,
let $r_i \leq r$ and $y_i \in [\theta]^{< \aleph_0}$ be such that
$r_i \vDash ``x_i = y_i"$.

$|{\Cal R}|<\theta$, therefore $(\exists \beta_0 < \kappa)
(\forall \beta \geq \beta_0)(\exists p_\beta \leq r)$ such that
$|\{i \in [\theta_{\beta},\theta_{\beta}^+) : p_\beta = r_i\}|=
\theta_{\beta}^+$.

Now, for every $\alpha < \kappa$, let $r_\alpha = p_{\beta_0 + \alpha}$ and
$$
T_\alpha=\{i \in [\theta_{\beta_0 + \alpha},\theta_{\beta_0 + \alpha}^+) : 
r_\alpha = r_i\}.
$$
Notice that for every $\alpha < \kappa$, 
$\{y_i : i \in T_\alpha\}$ are disjoint.  
Furthermore we may assume that 
$\{y_i : i \in \bigcup_{\alpha<\kappa}T_\alpha\}$
are disjoint (otherwise, by induction on $\alpha < \kappa$, pass to a subset
of $T_\alpha$ of cardinality $\theta_{\alpha}^+$). 
If $i \notin \bigcup_{\alpha<\kappa}T_\alpha$, 
then redefine $y_i=\emptyset$.  We now have $Y=\{y_i : i \in \theta \}$
disjoint finite subsets of $\theta$, and $\{ T_\alpha : \alpha < \kappa \}$
a $Y$-sequence.

Let $B$ be an ${\Cal R}$-name for the set 
$\{\alpha < \kappa : r_\alpha \in G \}$, where $G$ is a name for the
generic filter.  ${\Cal R}$ is $\kappa$.c.c., therefore there exists
$p \leq r$ such that $p \vDash ``|B|=\kappa"$.  Finally it is not hard 
to check that
$$
p \vDash ``\{\{\alpha \} : \alpha \in B \} \,\,\text {is an
antichain in}\,\, {\Cal Q}(Y,T)".
\qed
$$
\enddemo

\proclaim {Corollary} 
If $({\Cal A}, {\Cal B})$ is a nice $(\theta ,\theta)$ pair, then
for every $X= \{x_i : i \in \theta \}$ disjoint finite subsets of $\theta$,
and every $S=\{S_{\alpha} : \alpha < \kappa \}$ an $X$-sequence, 
${\Cal Q}(X,S)$ is not dangerous for $({\Cal A}, {\Cal B})$.
\endproclaim

\proclaim {Lemma 3} (See Lemma 8 [K1].)
Let $({\Cal A}, {\Cal B})$ be a nice $(\theta, \theta)$ pair.
Let $\gamma$ be a limit ordinal and ${\Cal P}_{\gamma}$ a finite support
iteration of c.c.c.\ partial orders.  If ${\Cal P}_{\gamma}$ is dangerous
for $({\Cal A}, {\Cal B})$,
then there exists $\alpha < \gamma$ such that ${\Cal P}_{\alpha}$ is dangerous
for $({\Cal A}, {\Cal B})$. 
\endproclaim

\demo {Proof}
Assume that ${\Cal P}_{\gamma}$ is dangerous for $({\Cal A}, {\Cal B})$.  
There are two cases.

$cof(\gamma) \neq \kappa$:  Let $X= \{x_i : i \in \theta \}\in 
\text{V}[G_\gamma]$, 
disjoint finite subsets of $\theta$, be a witness for the failure of
niceness,
where $G_\gamma$ is a ${\Cal P}_{\gamma}$-generic filter.  Now, there are three
subcases: $cof(\gamma) > \theta$, $\kappa < cof(\gamma) < \theta$, and
$cof(\gamma) < \kappa$. It is not hard to see that in each of these subcases
there exists $\alpha < \gamma$
and there is $A \in [\theta]^\theta$
such that $\{x_i : i \in A \}\in \text{V}[G_\alpha]$, contradicting
niceness in $\text{V}[G_\alpha]$.

$cof(\gamma) = \kappa$: Let $X$ be as in the previous case.
Let $c=\cup\{a(x_i) : i \in \theta \}\subset \omega$.
Let $\alpha < \gamma$ be such that $c \in \text{V}[G_\alpha]$.
(Such an $\alpha$ exists since $\kappa > \omega$.)
In $\text{V}[G_\alpha]$, let $Y$ be a disjoint family of finite subsets
of $\theta$ which is maximal with respect to the property
$$
(\forall y \in Y)\, [a(y) \subset c \wedge b(y) \cap c = \emptyset].
$$
$Y$ remains maximal in $\text{V}[G_\gamma]$ and therefore
must have cardinality $\theta$, which contradicts niceness in 
$\text{V}[G_\alpha]$.
\qed
\enddemo

\proclaim {Theorem {\rm (Kunen)}}
It is consisitent to have $m$ singular with $cof(m)=\omega_1$.
\endproclaim

\demo {Proof}
Let $\theta$ be singular with $cof(\theta)=\kappa=\omega_1$ and force with 
$\Cal P$ (see \S0) to start with 
$({\Cal A}, {\Cal B})$, a nice, disjoint, locally split $(\theta,\theta)$ pair. 
Let us now iterate c.c.c\. partial orders of size $< \theta$ in the 
following way.  Assume that the part of the iteration that has been defined 
thus far is non-dangerous.  Assume
that the next partial order on the list 
(of all c.c.c\. partial orders of size $< \theta$)
is $\Cal R$, but $\Cal R$ is dangerous
(otherwise just force with $\Cal R$).
Then by Lemma 2, there are $\Cal Y$ and $\Cal T$ 
such that ${\Cal R} \times {\Cal Q(Y,T)}$ is not c.c.c.  
In addition, by Lemma 1, ${\Cal Q(Y,T)}$ is c.c.c\. and by the  
corollary, ${\Cal Q(Y,T)}$ is non-dangerous for 
$({\Cal A}, {\Cal B})$.  So instead of forcing with 
$\Cal R$, let us force with ${\Cal Q(Y,T)}$ to add an uncountable antichain
to $\Cal R$.
By Lemma 3,the iteration of non-dangerous c.c.c\. partial orders
is a non-dangerous c.c.c\. partial order.  Therefore, in the extension, 
$({\Cal A}, {\Cal B})$ remains nice and $m=\theta$. 
\qed
\enddemo

In the general case ($\kappa$ is any regular uncountable cardinal), 
all we know is that ${\Cal Q(Y,T)}$ is $\kappa$.c.c.\ and not necessarily c.c.c.
So if $\kappa > \omega_1$, and we perform the iteration as in the proof 
of the theorem, then cardinals below $\kappa$ may be collapsed, and we may 
end up with a model for Kunen's result, in which $\kappa = \omega_1$.

\subhead \S2 Beyond niceness \endsubhead

Let us define a condition which implies niceness, and which is, in the presence
of MA$_\kappa$, equivalent to niceness.

\definition {Definition 2} 
Let $({\Cal A}, {\Cal B})$ be a nice $(\theta,\theta)$ pair.  We say 
that $(\ast)$ holds for $ ({\Cal A}, {\Cal B})$ 
if there is no c.c.c.\ partial order of cardinality $< \theta$
which is dangerous for $({\Cal A}, {\Cal B})$.
\enddefinition

The following Lemma shows that $(\ast)$ is in fact a statement concerning the
existence of certain dangerous
c.c.c.\ suborders of the various ${\Cal P}(Y,T)$'s.

\proclaim {Lemma 4}
Assume that $({\Cal A}, {\Cal B})$ is a nice $(\theta ,\theta)$ pair. Then 
$(\ast)$ fails for $({\Cal A}, {\Cal B})$ 
if and only if there are 
$Y= \{y_i : i \in \theta \}$ disjoint finite subsets of $\theta$,
and $T=\{T_{\alpha} : \alpha < \kappa \}$ a $Y$-sequence,
and there is ${\Cal P}' \in [{\Cal P}(Y,T)]^{\kappa}$
such that ${\Cal P}'$, equipt with the ordering of ${\Cal P}(Y,T)$,
is c.c.c.\  and closed under subsets. 
\endproclaim

\demo{Proof}
($\Leftarrow$): ${\Cal P}' \times {\Cal Q}(Y,T)$ is not $\kappa$.c.c.\ 
and hence, by Lemma 1, ${\Cal P}'$ is dangerous for $({\Cal A}, {\Cal B})$.
\newline
($\Rightarrow$): Let ${\Cal R}$ be a c.c.c.\ dangerous partial order of
cardinality $< \theta$.   
By the proof of Lemma 2, there are 
$Y= \{y_i : i \in \theta \}$ disjoint finite subsets of $\theta$,
and $T=\{T_{\alpha} : \alpha < \kappa \}$ a $Y$-sequence,
and there are $B$ an ${\Cal R}$-name for a subset of $\kappa$, and 
$p \in {\Cal R}$ such that 
$$
p \vDash ``|B|=\kappa \text{ and} \,\, 
\{\{\alpha \} : \alpha \in B \} \,\,\text{ is an
antichain in }\,\, {\Cal Q}(Y,T)".
$$
Let ${\Cal B(R)}$ be the boolean completion of ${\Cal R}$.
For every $F \in [\kappa ]^{<\aleph_0}$ let $r(F)=[\![ F \subset B ]\!] \cdot p$,
where $[\![ F \subset B ]\!]$ is the boolean value of 
$``F \subset B"$ in ${\Cal B(R)}$.
Now define 
$$
{\Cal P}' = \{ F \in [\kappa ]^{<\aleph_0} : r(F)>0 \}.
$$
${\Cal P}' \subset {\Cal P}(Y,T)$, and $|{\Cal P}'|=\kappa$
(because $p \vDash ``|B|=\kappa"$.)

Finally, assume that $F \perp_{{\Cal P}'} K$. Then $r(F \cup K)=0$, and 
therefore $r(F) \cdot r(K) =0$.  This proves that ${\Cal P}'$ is 
c.c.c.\  since
${\Cal B(R)}$ is c.c.c.
\qed
\enddemo

Similar to the proof of Fact 0, and using Lemma 4, one can now show that the 
nice disjoint locally split pair $({\Cal A},{\Cal B})$, that was added using 
${\Cal P}$ in \S 0, also satisfies the property $(\ast)$.  Let us check how 
well property $(\ast)$ is preserved through a c.c.c\. iteration of partial 
orders of size $<\theta$.  

We first show that $(\ast)$ is preserved through successor steps of the
iteration.  More precisely, if $(\ast)$ holds for $({\Cal A},{\Cal B})$ 
and ${\Cal R}$ is a c.c.c.\ partial
order with $|{\Cal R}|<\theta$, then 
$\vDash_{\Cal R} `` (\ast)$ holds for $({\Cal A},{\Cal B}) "$.
To show this assume otherwise.  By Lemma 4, there is $r \in {\Cal R}$
and $\pi$ an ${\Cal R}$-name such that
$$
\align
r\vDash_{\Cal R} ``&\text{there are}\,\, Y=\{y_i : i \in \theta \}
\,\,\text{ disjoint finite subsets of}\,\, \theta,\,\, \text {and}\,\, T, \,\,
\text{a}\,\, Y-\text{sequence 
with} \\
&\pi \in [{\Cal P}(Y,T)]^{\kappa},\,\,\text{and} \,\,\pi \,\, \text{is c.c.c.\  
and closed under subsets}".
\endalign
$$
In particular ${\Cal R}\ast \pi$ is a dangerous c.c.c.\ partial order.
But ${\Cal R}\ast \pi$ has a dense subset of cardinality $< \theta$,
namely $\{(r,F) : r\in {\Cal R}, F \in [\kappa]^{<\aleph_0}\,\, and \,\,
r\vDash'' F \in \pi "\}$. 
Therefore $(\ast)$ fails for $({\Cal A}, {\Cal B})$, which is a contradiction.

As for the preservation of $(\ast)$ at limits, we can only show the cases where
$cof(\gamma) \neq \kappa$ (see Lemma 3).  We first remark that in these cases
the following holds:
if $B \in [\kappa]^\kappa$, $Y= \{y_i : i \in \theta \}$
disjoint finite subsets of $\theta$,
and $T=\{T_{\alpha} : \alpha < \kappa \}$ a $Y$-sequence, are all in
$\text{V}[G_\gamma]$, then there exists
$\beta < \gamma$ such that in $\text{V}[G_\beta]$, there exists
$A \in [B]^\kappa$, and for every $\alpha \in A$, 
$T_{\alpha}' \in [T_{\alpha}]^{\theta_\alpha^+}$ are such that 
$\{y_i : i \in \bigcup_{\alpha\in A} T_{\alpha}' \} \in \text{V}[G_\beta]$.

The proof that this remark implies that $(\ast)$ is preserved proceeds as
follows.  Assume that $(\ast)$ fails in $\text{V}[G_\gamma]$.  
By Lemma 4, let $Y, T$ be given and
${\Cal P}' \in [{\Cal P}(Y,T)]^{\kappa}$
such that ${\Cal P}'$ is c.c.c.\  and closed under subsets.
Let $B=\{\alpha \in \kappa : \{\alpha\} \in {\Cal P}'\}$.
Now let $\beta < \gamma$ and $A \in [B]^\kappa$ as discussed above.
In $\text{V}[G_\beta]$ define the partial order
$$
\align
{\Cal P}'' = \{&F \in [A]^{< \aleph_0} : 
(\forall \alpha \neq \beta \in F) (\forall i \in T_{\alpha}')
(\forall j \in T_{\beta}') \\
&[a(y_i) \cap b(y_j) = \emptyset \,\wedge \, 
b(y_i) \cap a(y_j) = \emptyset ]\}.
\endalign
$$
In $\text{V}[G_\gamma]$, consider 
${\Cal R}' = \{F \in [A]^{< \aleph_0} : F \in {\Cal P}'\}$.
${\Cal R}' \subset {\Cal P}''$ and ${\Cal R}'$ is a c.c.c.\ partial order 
of cardinality $\kappa$.  
Finally, $\text{V}[G_\gamma]$ is a forcing extension of
$\text{V}[G_\beta]$, therefore in $\text{V}[G_\beta]$ we can define 
the partial order 
${\Cal R}''=\{ F \in {\Cal P}'' : [\![F \in {\Cal P}']\!] \cdot 
[\![{\Cal P}'\, \text{is c.c.c.\ } ]\!] > 0 \}$.
$|{\Cal R}''|=\kappa$ because ${\Cal R}' \subset {\Cal R}''$.
${\Cal R}''$ is c.c.c.\ because $\text{V}[G_\gamma]$ is a c.c.c.\ extension of
$\text{V}[G_\beta]$.  
Therefore, $\text{V}[G_\beta] \models ``{\Cal R}''\,\,\text{is c.c.c.\ and 
dangerous for}
\,\, ({\Cal A},{\Cal B})" $, and hence $(\ast)$ fails in $\text{V}[G_\beta]$.

Finally, let us look again at the case where $cof(\theta)=\kappa = \omega_1$.  
Let $({\Cal A}, {\Cal B})$ be a nice $(\theta ,\theta)$ pair.  
We claim that there exists a c.c.c.\ partial order ${\Cal Q}$ such that
$$
\vDash_{\Cal Q} ``(\ast)\,\,\text{holds for}\,\,({\Cal A}, {\Cal B})".
$$
${\Cal Q}$ is simply the finite support iteration of all partial orders of the
form $({\Cal Q}(X,S))^\omega$ (product with finite support).  ${\Cal Q}$ is 
c.c.c.\ and, by the Corollary and Lemma 3, it preserves the niceness of 
$({\Cal A}, {\Cal B})$.  In the extension, all partial orders of the form
${\Cal Q}(X,S)$ are $\sigma$-centered and therefore, by Lemma 2,  
$(\ast)$ holds for $({\Cal A}, {\Cal B})$.

This discussion suggests an alternative way of viewing Kunen's result.  
Start with a disjoint, locally split $(\theta,\theta)$ pair 
$({\Cal A},{\Cal B})$ for which $(\ast)$ holds (add such a pair using 
${\Cal P}$ which was defined in \S0).  Then iterate all c.c.c.\  partial orders 
of cardinality $< \theta$.  By the remarks above, 
$(\ast)$ may first fail only at limits of cofinality $\omega_1$. In this case,
by Lemma 3,
$({\Cal A}, {\Cal B})$ is still nice so   
first force with ${\Cal Q}$, as defined above, to get an extension in which
$(\ast)$ holds for $({\Cal A}, {\Cal B})$, and then force 
with the next c.c.c.\ partial order of cardinality $< \theta$ on the list.

\remark {Concluding remarks}
Given a disjoint $(\theta,\theta)$ pair $({\Cal A}, {\Cal B})$
we discussed two properties: \newline
{\tenrm(1)}\,$(\ast)$ holds for $({\Cal A},{\Cal B})$.\newline 
{\tenrm(2)}\,$({\Cal A},{\Cal B})$ is nice.

{\tenrm(1)} implies {\tenrm(2)}.   
{\tenrm(1)} is preserved while forcing with c.c.c.\  partial orders of 
cardinality less than $\theta$ (we do not have a similar result for 
{\tenrm(2)}). Both {\tenrm(1)} and {\tenrm(2)} are 
preserved at limit stages of c.c.c.\  iterations, of cofinality $\neq \kappa$.
{\tenrm(2)} is also preserved at limit stages of cofinality $= \kappa$ (we do 
not have a similar result for {\tenrm(1)}). 

Roughly speaking, it seems desirable to have uncountable c.c.c.\ 
suborders of the various ${\Cal Q}(X,S)$'s. This would enable us
to kill the dangerous partial orders as they come along,
or alternatively, force $(\ast)$ to hold and cosequently kill all the
dangerous partial orders at once.

It should be mentioned that if one starts with a nice $(\theta ,\theta)$ pair
$({\Cal A},{\Cal B})$, and then tries to preserve the niceness along
a c.c.c.\  iteration of all partial orders of size less than $\theta$, then 
the only  difficulty lies in getting MA$_\kappa$ to hold.  This is true because
if MA$_\kappa$ holds, then
$({\Cal A},{\Cal B})$ is nice if and only if 
$(\ast)$ holds for $({\Cal A},{\Cal B})$, 
and therefore one is free to force with the next partial order on the list.

On the otherhand if $\lambda < \kappa$, then MA$_\lambda$ could be forced 
without destroying niceness because c.c.c.\ partail orders of size $\lambda$ 
are not dangerous.  So if $\kappa > \omega_1$, then one can force 
MA$_{\aleph_1}$ and preserve niceness.  In this stage all c.c.c.\  partail 
orders are c.c.c.-productive and if there are still dangerous ones, they can 
not be killed and the iteration is stuck.   

\endremark


\Refs 
\ref \key {K}  \by K. Kunen  \book Set Theory \publ North-Holland
\yr1980
\endref
\ref \key {K1}  \by K. Kunen \yr1988 \paper Where MA First Fails 
\vol 53 \jour J. Symbolic Logic
\endref
\ref \key {S} \by S. Shelah \paper personal conversation  
\endref
\endRefs
\enddocument